\documentclass[reqno]{amsart}
\usepackage{epsfig,amsmath,amsfonts,latexsym}
\usepackage{amsthm}
\usepackage[abs]{overpic}
\usepackage[usenames,dvipsnames]{xcolor}
\usepackage{palatino}
\usepackage[hyperfootnotes=false]{hyperref}
\usepackage{dsfont}
\usepackage{comment}
\usepackage{accents}

\usepackage{amssymb}
\usepackage{braket}
\usepackage{adjustbox}
\usepackage{changepage}

\usepackage[utf8]{inputenc} 
\usepackage[T1]{fontenc} 

\usepackage{graphicx}
\usepackage{mathtools}
\usepackage{pifont}
\usepackage{tikz-cd}
\hypersetup{
  colorlinks,
  citecolor=Blue,
  linkcolor=Black,
  urlcolor=arsenic}
\usetikzlibrary{calc, 3d}

\definecolor{pakistangreen}{rgb}{0.0, 0.4, 0.0}
\definecolor{defnn}{RGB}{5, 66, 114} 
\definecolor{nicered}{HTML}{BA060F}

\theoremstyle{plain}

\newtheorem{mainthm}{Theorem}

\newtheorem{dummy}{anything}[section]
\newtheorem{thm}{Theorem}
\newtheorem{lemma}[dummy]{Lemma}

\newtheorem{prop}[dummy]{Proposition}

\newtheorem*{conj*}{Conjecture}
\theoremstyle{definition}
\newtheorem{definition}[dummy]{Definition}

\newtheorem{remark}[dummy]{Remark}

\newtheorem{assumption}{Assumption}
\theoremstyle{remark}

\newcommand{\M}{\mathcal{M}}

\newcommand{\norm}[1]{\left \lVert #1 \right \rVert}

\definecolor{pakistangreen}{rgb}{0.0, 0.4, 0.0}
\definecolor{myrtle}{rgb}{0.13, 0.26, 0.12} 
\definecolor{lincolngreen}{rgb}{0.11, 0.35, 0.02}
\definecolor{nicegreen}{HTML}{317041}
\definecolor{mintcream}{rgb}{0.96, 1.0, 0.98}
\definecolor{defn}{RGB}{5, 66, 114} 
\definecolor{lightimperial}{RGB}{119, 141, 179}
\definecolor{mauvelous}{rgb}{0.94, 0.6, 0.67}
\definecolor{cornflowerblue}{rgb}{0.39, 0.58, 0.93}
\definecolor{darkturquoise}{rgb}{0.0, 0.81, 0.82}
\definecolor{zapata}{HTML}{DEE9F5}
\definecolor{ao}{rgb}{0.16, 0.32, 0.75}
\definecolor{seacolor}{HTML}{6868D1}
\definecolor{nicered}{HTML}{BA060F}
\definecolor{amaranth}{rgb}{0.9, 0.17, 0.31}
\definecolor{nicepink}{HTML}{EB6EE6}
\definecolor{darkpastelpurple}{rgb}{0.59, 0.44, 0.84}
\definecolor{lavenderpurple}{rgb}{0.59, 0.48, 0.71}
\definecolor{brilliantlavender}{rgb}{0.96, 0.73, 1.0}
\definecolor{palatinatepurple}{rgb}{0.41, 0.16, 0.38}
\definecolor{purpleheart}{rgb}{0.41, 0.21, 0.61}
\definecolor{purple(x11)}{rgb}{0.63, 0.36, 0.94}
\definecolor{purpletaupe}{rgb}{0.31, 0.25, 0.3}
\definecolor{nicepurple}{HTML}{532E7B}
\definecolor{lavenderweb}{rgb}{0.9, 0.9, 0.98}
\definecolor{amber}{rgb}{1.0, 0.75, 0.0}
\definecolor{bronze}{rgb}{0.8, 0.5, 0.2}
\definecolor{deepcarrotorange}{rgb}{0.91, 0.41, 0.17}
\definecolor{coolorange}{rgb}{1.0, 0.49, 0.0}
\definecolor{liver}{rgb}{0.33, 0.29, 0.31}
\definecolor{dimgray}{rgb}{0.41, 0.41, 0.41}
\definecolor{cadet}{rgb}{0.33, 0.41, 0.47}
\definecolor{outerspace}{rgb}{0.25, 0.29, 0.3}
\definecolor{isabelline}{rgb}{0.96, 0.94, 0.93}
\definecolor{lightgray}{rgb}{0.83, 0.83, 0.83}
\definecolor{platinum}{rgb}{0.9, 0.89, 0.89}
\definecolor{moon}{HTML}{B3B3B3}
\definecolor{darkmoon}{HTML}{969696}
\definecolor{arsenic}{rgb}{0.23, 0.27, 0.29}
\definecolor{isabelline}{rgb}{0.96, 0.94, 0.93}
\definecolor{linen}{rgb}{0.98, 0.94, 0.9}
\definecolor{snow}{rgb}{1.0, 0.98, 0.98}

\renewcommand{\indent}{\hspace{1em}}

\title{Bi-normal trajectories in the Circular Restricted Three-Body Problem}
\author{Agustin Moreno, Arthur Limoge}

\address[A.\ Moreno, A.\ Limoge]{Heidelberg University, DE}

\email{agustin.moreno2191@gmail.com, arthur.limoge@outlook.com}

\date{}

\begin{document}

\begin{abstract}
In this note, we study existence of infinitely many trajectories \emph{bi-normal} (i.e.\ normal at initial and final times) to the $xz$-plane in the Spatial Circular Restricted Three-Body problem, in the convexity range and near the primaries, under the assumption of the twist condition as defined by Moreno--van-Koert in \cite{Mor01}. Modulo our assumptions, this is an expected application of the relative Poincaré--Birkhoff theorem for Lagrangians in Liouville domains, proven by the authors in \cite{ML24}.

\end{abstract}

\maketitle

\tableofcontents

\section{Introduction}

Let $\mathbb{R}^3$ denote three-dimensional space with coordinates $q_1, q_2, q_3$, and consider three bodies, the Earth ($E$), the Moon ($M$), and a satellite ($S$), under the mutual influence of Newtonian gravity. We assume that the Earth and Moon move in circles about each other, so that we can choose a rotating frame and fix their positions in the plane. We also assume the satellite has mass $m_S = 0$. This system is called the (spatial) circular restricted $3$-body problem, or SCR3BP for short.

\begin{definition}[\emph{Bi-normal trajectories}]\label{q2symmetryDefn}
    A trajectory $x(t)=\big(q(t),\dot{q}(t)\big)$ is \emph{bi-normal} to the $xz$-plane if there exist times $t_0<t_1$ for which $x(t_0)$ and $x(t_1)$ are normal to the plane $\{q_2=0\}$, i.e.\ $$q_2(t_j)=\dot{q}_1(t_j)=\dot{q}_3(t_j)=0,$$ for $j=0,1.$
\end{definition}

\vspace{-1.2em}

\[
\begin{adjustbox}{scale=1.3}
\begin{tikzpicture}[scale=2]

    \fill[snow] (-2,0,0) -- (2,0,0) -- (2,0,3) -- (-2,0,3) -- cycle;
    \draw (-2,0,0) -- (2,0,0);
    \draw (2,0,0) -- (2,0,3);
    \draw (2,0,3) -- (-2,0,3);
    \draw (-2,0,3) -- (-2,0,0);

    \begin{scope}[shift={(-0.6,-0.05,1)}] 
        \def\E{0.30} 
        \fill[blue!70!black!70] (0,0) circle (\E);
        \draw[thick] (0,0) circle (\E);
        \draw[very thin, ball color=blue!70!black!40, fill opacity=0.2] (0,0) circle (\E);
        \begin{scope}[rotate=-11]
            \clip (0,0) circle (\E);
            \fill[white] (0,\E) ellipse ({0.6*\E} and {0.15*\E});
            \fill[white] (0,-\E) ellipse ({0.8*\E} and {0.08*\E});
            \fill[green!70!black!60,rotate=-30] (160:1.1*\E) ellipse ({0.2*\E} and {0.8*\E});
            \fill[green!70!black!60,rotate=40] (-10:1.14*\E) ellipse ({0.2*\E} and {0.9*\E});
            \fill[green!60!black!60,very thick,rotate=-20] 
                (230:0.86*\E) ellipse ({0.25*\E} and {0.18*\E});
        \end{scope}
    \end{scope}

    \begin{scope}[shift={(1,-0.05,1)}] 
        \def\M{0.18} 
        \fill[gray!70!white] (0,0) circle (\M);
        \draw[thick] (0,0) circle (\M); 
        \draw[very thin, ball color=gray!70!white, fill opacity=0.2] (0,0) circle (\M);
        \fill[gray!50!white] (0.05,0.03) circle (0.025);
        \fill[gray!50!white] (-0.04,-0.05) circle (0.03);
        \fill[gray!60!white] (-0.02,0.05) circle (0.02); 
        \fill[gray!60!white] (0.04,-0.02) circle (0.015); 
        \fill[gray!40!white] (-0.05,-0.05) circle (0.04); 
        \fill[gray!40!white] (0.02,0.05) circle (0.035);
        \fill[gray!40!white] (0.01,-0.06) circle (0.025); 
        \fill[gray!60!white] (-0.06,0.02) circle (0.018);
    \end{scope}

    \node at (0.2,0.5,1) {\includegraphics[width=0.1\textwidth]{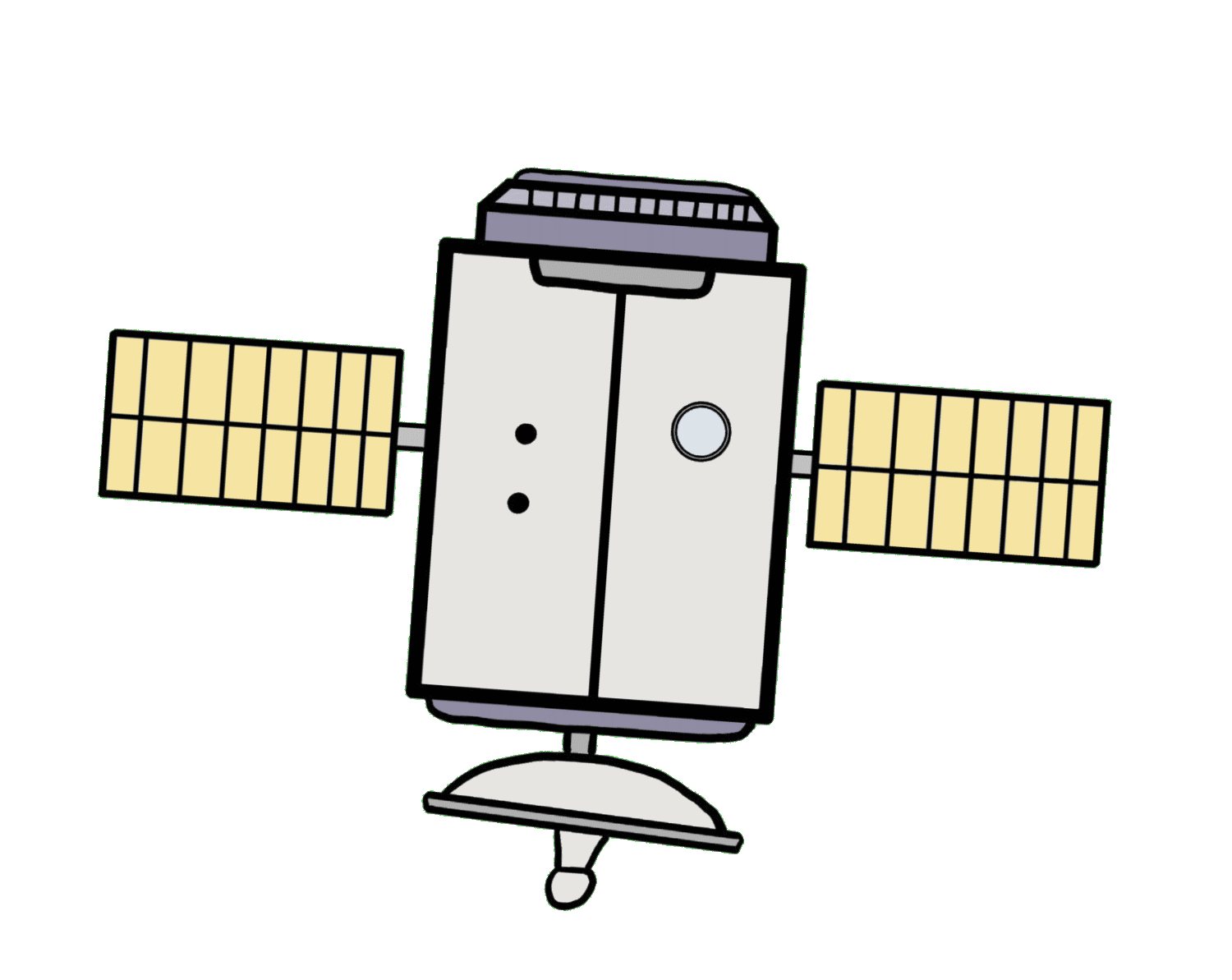}};

    \draw[thick, nicered,->] (0.2,0.5,1) -- (0.2,0.5,2); 

    \node at (-1.95, -1.05, 0) {\small{Earth-Moon (ecliptic) plane}};
    \draw[dashed] (-1.7,0,0.2) -- (-1.4,0,0.2) node[anchor=north] {\small $q_1$};
    \draw[dashed] (-1.7,0,0.2) -- (-1.7,0,0.5) node[anchor=east] {\small $q_2$};
    \draw[->, thick] (-1.7,0,0.2) -- (-1.4,0,0.2); 
    \draw[->, thick] (-1.7,0,0.2) -- (-1.7,0,0.5);
    \draw[dashed] (-1.7,0,0.2) -- (-1.7,0.3,0.2) node[anchor=east] {\small $q_3$}; 
\draw[->, thick] (-1.7,0,0.2) -- (-1.7,0.3,0.2);  
\end{tikzpicture}
\end{adjustbox}
\]

In the figure, the $q_1$-axis corresponds to the Earth-Moon axis. We have drawn the position of the satellite at time $t = t_0$, normal to the $xz$-plane (or $q_1q_3$-plane, in our coordinates).

\medskip

There are $5$ critical points of the Hamiltonian describing the SCR3BP, the \emph{Lagrangians} $L_1,\dots,L_5$, ordered by increasing energy. We call \emph{low energy range} the range of energies below the energy of $L_1$, or \emph{slightly} above it (so that the main results of \cite{AFvKP,CJK} apply). For this range of energies, the \emph{Hill region} (i.e.\ the region the satellite $S$ can move in for its given energy) has either 3 or 2 connected components respectively, one of which is unbounded. We will refer to the dynamics in the bounded components (containing the Earth and Moon) as the dynamics \emph{near the primaries}. 

\medskip

\textbf{Main result.} The goal of this note is to explain the heuristics behind the following speculative result.

\begin{mainthm} \textbf{(speculative)} \label{MainThm}
    Assuming the twist condition of Moreno--van-Koert (Assumption \ref{TwistCondition}) holds in the SCR3BP, then in the convexity range, and near the primaries, there exist infinitely many trajectories which are bi-normal to the $xz$-plane.
\end{mainthm}

Modulo some technicalities which we explain below, this is an expected application of (a modified version of) the relative Poincaré--Birkhoff theorem from \cite{ML24}, proving the existence of infinitely many Hamiltonian chords, given a Lagrangian with infinite-dimensional wrapped Floer cohomology in a Liouville domain. An analogous statement for consecutive spatial collision orbits was already obtained in \cite{ML24} via the same Floer-theoretical methods, as adaptations of the arguments from Ginzburg's proof of the Conley conjecture \cite{Conley_conjecture}. 

\indent The 'convexity range' mentioned in the theorem is a subset of the low energy range. It is the range of pairs $(\mu, c)$ such that the Levi-Civita regularization of the planar CR3BP is convex, where $\mu = m_M/(m_M+m_E)$ is the Earth-Moon normalized mass-ratio, and $c$ the Jacobi constant (the energy). It is an \textit{a priori} non-perturbative set of values; see \cite{Albers2012}.

\smallskip

\textbf{A twist condition.} To state the twist condition of Moreno--van-Koert, let us recall some theory about the SCR3BP. For energy below that of $L_1$, near the primaries and after regularization of collisions, the SCR3BP flow can be viewed as a Reeb flow on a fibrewise star-shaped domain in $T^\star\mathbb{S}^3$, by \cite{CJK} (or $T^\star\mathbb{S}^3\natural T^\star\mathbb{S}^3$ for energy slightly above $H(L_1)$). Moreover, by \cite{Mor00}, one can find an open book decomposition of the regularized energy level set which is adapted to the flow. Each page of the open book is then a global hypersurface of section symplectomorphic to a \emph{degenerate} Liouville domain $(W, \omega=\mathrm{d}\lambda)$. Here, $\lambda$ is the restriction to the page of the ambient contact form, and therefore the symplectic form becomes degenerate at the boundary (this is what we mean by a degenerate Liouville domain).

Since $W$ is a global hypersurface of section, we can consider the associated Poincar\'e return map $$\tau : \mbox{int}(W) \to \mbox{int}(W),$$ mapping each point to its first return point under the flow. This is a Hamiltonian diffeomorphism in the interior under convexity assumptions, which hold for the convexity range (see \cite{Mor00}). This map was also shown to extend \emph{smoothly} to the boundary in \cite{Mor00} under the same convexity assumptions.

\indent There is also an alternative description of the setup which can be obtained by conjugating by a map which is smooth in the interior of the page, but only continuous at the boundary. The result is that the $2$-form $\omega$ can be pulled back by this map to actually become symplectic along the boundary, but the price to pay is that the return map only extends continuously to the boundary. To summarize, one can realistically choose either of the following equivalent setups:

\begin{itemize}
    \item[(a)] The return map $\tau$ extends smoothly to the boundary, but $\omega$ degenerates there; or
    \item[(b)] The return map $\tau$ extends only continuously to the boundary, but $\omega$ is actually non-degenerate everywhere.
\end{itemize}

The above dichotomy is explained in detail in Moreno's recent book, see Section 6 of \cite{moreno2024symplecticgeometrythreebodyproblem}. In what follows, as is done in \cite{Mor01}, we will make the simplifying (although unrealistic) assumption that $\tau$ extends \emph{smoothly} to the boundary, while $\omega$ is also symplectic everywhere; see Remark \ref{rk:difficulties} for more discussion on this issue. In other words, we will assume the above technicalities away, and work with Liouville domains and smooth maps. Moreover, we assume the following:

\begin{assumption}[Twist condition]\label{TwistCondition}
    The return map $\tau : W \to W$ is generated by a $\mathcal{C}^2$ Hamiltonian $H_t : W \to \mathbb{R}$, whose Hamiltonian vector field $X_{H_t}$ satisfies

    \vspace{-0.6em}

    \begin{equation}
        X_{H_t}|_{\partial W} = h_t \mathcal{R}_\alpha,
    \end{equation}

    \smallskip

    where $h_t > 0$ is a smooth function on $\partial W$, and $\mathcal{R}_\alpha$ is the Reeb vector field on $(\partial W, \alpha := \lambda|_{\partial W})$.
\end{assumption}

\begin{remark}[\textbf{Difficulties}]\label{rk:difficulties}
    This twist condition is a generalization to arbitrary dimension of Poincaré's original twist condition \cite{PLT}, formulated for the planar problem. For the spatial problem, however, it is quite a strong assumption, and at this point, though we know of Hamiltonians that generate $\tau$ and even extend to the boundary if we choose setup (a) above (see \cite{Mor01}), we do not know of any satisfying a twist condition.
    
    \indent As explained above, there is still the difficulty in that the symplectic form on $W$ degenerates at the boundary. While there is an obvious version of the twist condition for degenerate Liouville domains (i.e.\ simply consider maps which correspond to twist maps under the continuous conjugation), proving a fixed-point theorem in this context is a problem which has not yet been successfully addressed, as in setup (b) the Hamiltonian vector field explodes at the boundary (in the direction of the Reeb vector field; so that one would be looking at ``infinitely twisting'' maps). This phenomena is in fact present in the planar problem as well (see e.g.\ \cite{moreno2024symplecticgeometrythreebodyproblem} for more details on this difficulty). This implies that one needs an alternative, weaker notion of a twist condition (in particular, addressing the boundary degeneracy), but for which one can still obtain analogous results. This is subject of ongoing efforts, and is the reason why the above expected application is still speculative.
\end{remark}

\textbf{Idea for the proof.} The key point for the proof of our Main Theorem is that bi-normal trajectories can be seen as chords with ends in a Lagrangian, inside the global hypersurface of section $W$. This Lagrangian corresponds to the fixed-point locus of one of the symmetries of the SCR3BP (i.e.\ an anti-symplectic involution that preserves the Hamiltonian), and it can be interpreted as the co-normal bundle of a submanifold inside $W$. Therefore its wrapped Floer homology can be shown to be infinite-dimensional, and the main theorem from \cite{ML24} (suitably modified as explained above) should yield the above application.

\smallskip

\textbf{Acknowledgements.} For this work, the authors were supported by the Air Force Office of Scientific Research (AFOSR) under Award No.\ FA8655-24-1-7012, by the DFG under Germany's Excellence Strategy EXC 2181/1 - 390900948 (the Heidelberg STRUCTURES Excellence Cluster), and by the Sonderforschungsbereich TRR 191 Symplectic Structures in Geometry, Algebra and Dynamics, funded by the DFG (Projektnummer 281071066 – TRR 191). \

\section{Symmetries of the spatial problem}\label{SectionSymm}

After choosing a rotating frame in the Earth-Moon plane (= the \textit{ecliptic}) in order to fix their positions at $\vec{E}$ and $\vec{M}$ respectively, the Hamiltonian for the SCR3BP is:

\vspace{-1em}

\begin{align}
    H : T^\star(\mathbb{R}^3\backslash\{\vec{E},\vec{M}\}) &\longrightarrow \mathbb{R} \\
    (q,p) &\longmapsto \dfrac{1}{2}\norm{p}^2 - \dfrac{m_E}{\Vert q - \vec{E}\Vert} - \dfrac{m_M}{\Vert q - \vec{M}\Vert} + q_1p_2 - q_2p_1 \notag.
\end{align}

From this expression, we can write down three natural symmetries of the SCR3BP:

\begin{itemize}
    \item A symplectic involution $r : \mathbb{R}^6 \to \mathbb{R}^6 : (q_1, q_2, q_3, p_1, p_2, p_3) \mapsto (q_1, q_2, -q_3, p_1, p_2, -p_3)$. Its fixed point set is

    \vspace{-0.8em}

    \begin{equation}
        \text{Fix}(r) = \{(q,p) \in \mathbb{R}^6 \, \mid \, q_3 = p_3 = 0\}.
    \end{equation}

    \smallskip
    
    The symmetry $r$ is simply induced by reflection about the ecliptic $\{q_3=0\}\subset \mathbb R^3.$

    \medskip

    \item An anti-symplectic involution $\rho_1 : \mathbb{R}^6 \to \mathbb{R}^6 : (q_1, q_2, q_3, p_1, p_2, p_3) \mapsto (q_1, -q_2, -q_3, -p_1, p_2, p_3)$. Its fixed point set 

    \vspace{-0.8em}

    \begin{equation}
        \text{Fix}(\rho_1) = \{(q,p) \in \mathbb{R}^6 \, \mid \, q_2 = q_3 = p_1 = 0\},
    \end{equation}
    is Lagrangian in $T^\star(\mathbb{R}^3\backslash\{\vec{E},\vec{M}\})$.
    
    \smallskip

    \item An anti-symplectic involution $\rho_2 : \mathbb{R}^6 \to \mathbb{R}^6 : (q_1, q_2, q_3, p_1, p_2, p_3) \mapsto (q_1, -q_2, q_3, -p_1, p_2, -p_3)$. Its fixed point set

    \vspace{-0.8em}

    \begin{equation}\label{FixRho2Expr}
        \text{Fix}(\rho_2) = \{(q,p) \in \mathbb{R}^6 \, \mid \, q_2 = p_1 = p_3 = 0\},
    \end{equation}
    is also Lagrangian.
\end{itemize}

\bigskip

We will be interested in the last two symmetries, most particularly $\rho_2$. Indeed, observe that Definition \ref{q2symmetryDefn} can be rephrased by saying that a trajectory is bi-normal to the $xz$-plane iff it starts on $\text{Fix}(\rho_2)$ at time $t = t_0$, and comes back to it at $t = t_1$.






\bigskip

Now, to formulate this as a Floer-theoretical problem, we first need to regularize at collisions (indeed, notice that $H$ is singular as $q \to \vec{E}$ or $q \to \vec{M}$, causing energy hypersurfaces $H^{-1}(c)$ to be non-compact). To regularize, we employ Moser regularization, following §4 of \cite{Mor00}:

\vspace{-0.2em}

\begin{equation*}
    (q,p) \xrightarrow{\text{Switch map}} (p, -q) =: (x,y) \xrightarrow{\text{Inverse stereographic projection}} (\xi, \eta) \in T^\star \mathbb{S}^3\backslash\{N\}
\end{equation*}

\bigskip

Here, $T^\star\mathbb{S}^3$ is written in coordinates as:

\vspace{-0.4em}

\begin{equation}
    T^\star\mathbb{S}^3 := \left\{(\xi,\eta) \in T^\star\mathbb{R}^4 \, \mid \, \norm{\xi}^2 = 1, \braket{\xi,\eta} = 0\right\},
\end{equation}

\smallskip

and the North Pole is the point $N = (1, 0, 0, 0) \in \mathbb{S}^3$. The formula for inverse stereographic projection $T^\star\mathbb{R}^3 \to T^\star\mathbb{S}^3\backslash\{N\}$ is given by:

\vspace{-0.8em}

\begin{align}
    \xi_0 &= \dfrac{\norm{x}^2 - 1}{\norm{x}^2 + 1} \notag \\
    \xi_i &= \dfrac{2x_i}{\norm{x}^2 + 1} \, \, \text{ for i = 1,2,3} \\
    \eta_0 &= \braket{x,y} \notag \\
    \eta_i &= \dfrac{\norm{x}^2 + 1}{2}y_i - \braket{x,y} x_i \, \, \text{ for i = 1,2,3} \notag. 
\end{align}

\bigskip

With this, we can explicitly compute the regularizations $\widetilde{F}_1$ and $\widetilde{F}_2$ of $\text{Fix}(\rho_1)$ and $\text{Fix}(\rho_2)$:

\vspace{-0.8em}

\begin{align}
     \widetilde{F}_1 &= \left\{(\xi, \eta) \in T^\star\mathbb{S}^3 \, \mid \, \xi_1 = \eta_0 = \eta_2 = \eta_3 = 0\right\}, \\
    \widetilde{F}_2 &= \left\{(\xi, \eta) \in T^\star\mathbb{S}^3 \, \mid \, \xi_1 = \xi_3 = \eta_0 = \eta_2 = 0\right\},
\end{align}

\bigskip

which are submanifolds of $T^\star\mathbb{S}^3$. Moreover, after regularization, our Hamiltonian $H$ becomes

\vspace{-0.7em}

$$Q : T^\star\mathbb{S}^3 \longrightarrow \mathbb{R} : (\xi,\eta) \longmapsto \frac{1}{2}f(\xi,\eta)^2\norm{\eta}^2,$$

where $f$ is a positive function whose expression is given in §4 of \cite{Mor00}. Then, our regularized energy hypersurface corresponds to $\widetilde{\Sigma} := Q^{-1}(1/2) \cong \mathbb{S}^\star\mathbb{S}^3$, for energy below $H(L_1)$. And as we already mentioned in the Introduction, \cite{Mor00} constructs an open book on $\tilde{\Sigma}$ which is adapted to the SCR3BP dynamics. In particular, each of its pages is a global hypersurface of section. And (see §6 of \cite{Mor01}), this open book admits a particularly nice page given by
\vspace{-0.4em}
\begin{equation}
    W = \left\{(\xi, \eta) \in T^\star\mathbb{S}^3\, \mid \, \xi_3 = 0, \eta_3 \geq 0 \right\}.
\end{equation}

Similar remarks apply for energy slightly above $H(L_1)$, where the corresponding $W$ is a connected sum of two copies of the above. We will focus only on the former case in what follows, but all arguments carry over naturally to the latter. 

\medskip

The following was observed in \cite{Mor00}:

\begin{lemma}
    For energy below $H(L_1)$, $W$ is symplectomorphic to a fibre-wise star-shaped domain in $T^\star\mathbb{S}^2$, with its standard symplectic structure, and is therefore diffeomorphic to $\mathbb D^*\mathbb{S}^2$.
\end{lemma}

Here, a little care is needed near the boundary, as explained in Remark \ref{rk:difficulties} (i.e.\ in the above statement we consider setup (b) as explained above; we do the same below when describing the Lagrangians). The dynamics at the boundary $\partial W$ is the regularized dynamics of the planar CR3BP.

\medskip

\ding{226} Let us now consider the intersection of our regularized fixed point sets with the page $W$:

\vspace{-1.2em}

\begin{align*}
    \widetilde{F}_1\cap W &= \big\{(\xi, \eta) \in T^\star\mathbb{R}^4 \, \mid \, Q(\xi, \eta) = \frac{1}{2}, \norm{\xi}^2 = 1, \xi_1 = \xi_3 = \eta_0 = \eta_2 = \eta_3 = 0\big\} \\
    \widetilde{F}_2\cap W &= \big\{(\xi, \eta) \in T^\star\mathbb{R}^4 \, \mid \, Q(\xi, \eta) = \frac{1}{2}, \norm{\xi}^2 = 1, \xi_1 = \xi_3 = \eta_0 = \eta_2 = 0, \eta_3 \geq 0\big\}
\end{align*}

\smallskip

Note that in the second case, from the energy constraint $Q(\xi,\eta)=1/2$, the coordinate $\eta_3\geq 0$ is completely determined by the remaining ones. We rewrite:

\vspace{-0.8em}

\begin{align}
    \widetilde{F}_1 \cap W 
    &= \big\{(\xi, \eta) = (\xi_0, 0, \xi_2, 0, 0,  \eta_1, 0, 0) \, \mid \, \xi_0^2 + \xi_2^2 = 1, \eta_1^2 = 1/f^2(\xi_0,\xi_2) \big\}\\
    \tilde{F}_2\cap W &= \big\{(\xi, \eta) = (\xi_0, 0, \xi_2, 0, 0, \eta_1, 0, \eta_3) \, \mid \, \norm{\xi}^2 = 1, \eta_1^2+\eta_3^2 = 1/f^2(\xi_0,\xi_2), \eta_3 \geq 0\big\}
\end{align}

Meanwhile, the regularization of $\text{Fix}(r)$ gives $\{\xi_3 = 0, \eta_3 = 0\}$, which is the binding of the SCR3BP open book (or, in other words, the regularized planar problem).

\begin{definition}\label{defn:L2}
    We write $\Lambda_1 := \widetilde{F}_1\cap W$, and $L_2 := \widetilde{F}_2\cap W$. Observe that $\Lambda_1 = \partial L_2$. 
\end{definition}

\begin{prop}
    $L_2$ is an exact Lagrangian in $W$. Its boundary $\partial L_2 = \Lambda_1$ is Legendrian in $\partial W$.
\end{prop}

 \begin{adjustwidth}{0.5cm}{}
    \proof $W$ inherits its symplectic form from the embedding $W \hookrightarrow T^\star\mathbb{S}^3 \hookrightarrow T^\star\mathbb{R}^4$, hence

    \vspace{-0.8em}

    \begin{equation*}
        \omega|_W = \displaystyle\sum_{i=0}^2 \mathrm{d}\xi_i \wedge\mathrm{d}\eta_i\Big\vert_{W} = \mathrm{d}\left(- \displaystyle\sum_{i=0}^2 \eta_i \mathrm{d}\xi_i\right)\Bigg\vert_{W} = \mathrm{d}\lambda.
    \end{equation*}

$\lambda$ vanishes on $L_2$, and so the claim follows. \qed 
\end{adjustwidth}

\bigskip

In other words, $L_2$ is an exact Lagrangian with Legendrian boundary in $W$ (and it is also spin, because it is an orientable surface), making it admissible for wrapped Floer cohomology.

\section{Proof of Theorem \ref{MainThm}}\label{ProofSection}

Recall we want to prove Theorem \ref{MainThm}, namely that if the SCR3BP satisfies the twist condition, then there exist infinitely many trajectories bi-normal to the $xz$-plane (Definition \ref{q2symmetryDefn}) near the primaries, for low energy and in the convexity range. As we observed right after (\ref{FixRho2Expr}), such trajectories can be interpreted as Hamiltonian chords on the Lagrangian $L_2$.

\begin{definition}
    Let $H_t : W\to\mathbb{R}$ be a Hamiltonian. A \textit{Hamiltonian chord} $x$ of $H_t$ on a Lagrangian $L$ is a trajectory $x : [0,1] \to W$ of the Hamiltonian flow such that $x(0), x(1) \in L$.
\end{definition}

\smallskip

Hence, Theorem \ref{MainThm} should reduce to showing that there exist infinitely many such chords on $\text{int}(L_2)$. This should heuristically follow from a version of the following (where recall that we are assuming the boundary difficulties away):

\begin{thm}[\cite{ML24}]\label{MLPBThm}
    Let $(W, \omega = \mathrm{d}\lambda)$ be a connected Liouville domain, and $L \subset W$ be an exact, spin Lagrangian with Legendrian boundary. Let $\tau : W\to W$ be an exact symplectomorphism, and assume:

    \begin{itemize}
        \item $\tau$ satisfies the twist condition (Assumption \ref{TwistCondition}).

        \item $\dim HW^{*}(L) = \infty$, where $HW^{*}$ denotes \textbf{wrapped Floer cohomology}.

        \item if $\dim W > 2$, then: $c_1(TW) = 0$, and $(\partial W, \alpha)$ is strongly index-definite (see \cite{ML24}).

        \item $\tau$ admits finitely many interior periodic chords.
    \end{itemize}

    \smallskip

    Then, $\tau$ admits infinitely many interior Hamiltonian chords on $L$, of arbitrarily large order.
\end{thm}

\medskip

\begin{remark}
    We already know the strong index-definiteness assumption to hold in the SCR3BP, for energies $c < H(L_1) + \varepsilon$, and in the convexity range (see \cite{Mor01}). Potential weakenings of the twist condition, allowing one to also lift the strong-index definiteness assumption from Theorem \ref{MLPBThm} (and hence work outside of the convexity range) are discussed in \cite{Lim}.
    
    \indent Meanwhile, the fourth assumption is harmless: indeed, if it did not hold, then we would already have infinitely many trajectories (though maybe with bounded orders).
\end{remark}

The only thing we have left to do to prove Theorem \ref{MainThm} is to show that $\dim HW^{*}(L_2) = \infty$. First, observe that $L_2$ can be viewed as a conormal bundle in $\mathbb{D}^\star\mathbb{S}^2$ by:

\begin{lemma}[\cite{ASconormal}, Prop. 2.1]
    Let $M$ be a manifold, and $L$ a submanifold of $T^\star M$ on which the Liouville form $\lambda$ vanishes identically. Then the intersection of $L$ with the zero section of $T^\star M$ is a submanifold $R$. Furthermore, if $L$ is a closed subset of $T^\star M$, then $L = N^\star R$.
\end{lemma}

\smallskip

But we can actually be more precise. Indeed, by definition we have:

\vspace{-0.6em}

\begin{equation}\label{L2Cassandra}
    L_2 = \left\{(\tilde{\xi},\tilde{\eta}) = (\xi_0, 0, \xi_2, 0, \eta_1, 0) \, \mid \, \xi_0^2 + \xi_2^2 = 1, \eta_1^2 \leq 1/f^2\right\} \subset \mathbb{D}^\star\mathbb{S}^2
\end{equation}

\smallskip

And so we can observe:

\smallskip

\begin{lemma}
    $L_2 = N^\star R \cap W$, where $N^\star R \subset T^\star \mathbb{S}^2$ is the conormal bundle of the equator:

    \vspace{-0.6em}

    \begin{equation}\label{REquatorExpr}
        R := \left\{(\xi_0, 0, \xi_2) \in \mathbb{R}^3 \, \mid \, \xi_0^2 + \xi_2^2 = 1\right\} \subset \mathbb{S}^2.
    \end{equation}
\end{lemma}

\begin{adjustwidth}{0.5cm}{}
    \proof We directly compute:

        \vspace{-0.6em}

    \begin{equation*}
        N^\star R \cap W= \left\{(\xi_0, 0, \xi_2, 0, \eta_1, 0) \, \mid \, \xi_0^2 + \xi_2^2=1, \eta_1 \in \mathbb{R}) \right\} \cap W \cong L_2. 
    \end{equation*} \qed 
\end{adjustwidth}

\bigskip

To then compute wrapped Floer cohomology, we need to Liouville complete both $L_2$ and $W$, so that we get a Lagrangian $\widehat{L}_2 \cong N^\star\mathbb{S}^1$ inside $\widehat W \cong T^\star\mathbb{S}^2$. Then, to show that $HW^{*}(L_2)$ is infinite-dimensional, we invoke:

\begin{thm}[\cite{ASconormal}]\label{ConormalThm}
    Let $M$ be a closed, orientable manifold, and $L \subset T^\star M$ be a Lagrangian such that $L = N^\star R$ for some closed submanifold $R \subset M$. Then:

    \vspace{-0.4em}

    \begin{equation*}
        HW^{*}(L) \cong H^{*}(\mathcal{P}_RM)
    \end{equation*}

    \medskip

    where $\mathcal{P}_R M$ is the space of $\mathcal{C}^\infty$ paths in $M$ with endpoints in $R$; and $H^{*}$ denotes singular cohomology.
\end{thm}

\smallskip

allowing us to deduce:

\medskip

\underline{\textit{Proof of Theorem \ref{MainThm}.}} We have $HW^{*}(L_2) \cong H^{*}(\mathcal{P}_{\mathbb{S}^1}\mathbb{S}^2)$. This is infinite-dimensional, by standard methods of algebraic topology (see \cite{Lim}), allowing us to apply Theorem \ref{MLPBThm} and conclude. \qed 

\medskip

\begin{remark}[On trajectories binormal to the $x$-axis]
    By essence, Theorem \ref{MLPBThm} only counted chords in the \textit{interior} of our Liouville domain; which in our case, corresponded to trajectories bi-normal to the $xz$-plane. But $\partial L_2$ is also interesting, since it is the regularized fixed point set of the first symmetry. If we could prove a similar Poincaré-Birkhoff type theorem for chords on the boundary ($\partial W, \alpha)$, with respect to the Legendrian $\Lambda_1 = \partial L_2$, then we would obtain a result of the form:

    \begin{conj*}
        Assuming the twist condition or a variation thereof, there exist infinitely many trajectories \textbf{bi-normal to the $x$-axis} in the Spatial Circular Restricted Three-Body Problem, in the low energy range, and near the primaries. i.e there exist times $t_0 < t_1$ such that:

        \vspace{-0.8em}

        \begin{align*}
            q_2(t_j) = q_3(t_j) = \dot{q}_1(t_j) = 0
        \end{align*}

        for $j = 0,1$. In other words, the trajectory starts on the Earth-Moon axis, but with velocity pointing strictly outward, and comes back to satisfy the same condition after finite time.
    \end{conj*}

\[
\begin{adjustbox}{scale=1.5}
\begin{tikzpicture}[scale=2]

    \draw[thick, dashed] (-1.5,0,0) -- (1.5,0,0);

    \begin{scope}[shift={(-0.9,0,0)}] 
        \def\E{0.30} 
        \fill[blue!70!black!70] (0,0) circle (\E);
        \draw[thick] (0,0) circle (\E); 
        \draw[very thin, ball color=blue!70!black!40, fill opacity=0.2] (0,0) circle (\E);
        \begin{scope}[rotate=-11]
            \clip (0,0) circle (\E);
            \fill[white] (0,\E) ellipse ({0.6*\E} and {0.15*\E});
            \fill[white] (0,-\E) ellipse ({0.8*\E} and {0.08*\E});
            \fill[green!70!black!60,rotate=-30] (160:1.1*\E) ellipse ({0.2*\E} and {0.8*\E});
            \fill[green!70!black!60,rotate=40] (-10:1.14*\E) ellipse ({0.2*\E} and {0.9*\E});
            \fill[green!60!black!60,very thick,rotate=-20] 
                (230:0.86*\E) ellipse ({0.25*\E} and {0.18*\E});
        \end{scope}
    \end{scope}

    \begin{scope}[shift={(1.1,0,0)}]
        \def\M{0.18}
        \fill[gray!70!white] (0,0) circle (\M); 
        \draw[thick] (0,0) circle (\M); 
        \draw[very thin, ball color=gray!70!white, fill opacity=0.2] (0,0) circle (\M);

        \fill[gray!50!white] (0.05,0.03) circle (0.025);
        \fill[gray!50!white] (-0.04,-0.05) circle (0.03); 
        \fill[gray!60!white] (-0.02,0.05) circle (0.02); 
        \fill[gray!60!white] (0.04,-0.02) circle (0.015); 
        \fill[gray!40!white] (-0.05,-0.05) circle (0.04); 
        \fill[gray!40!white] (0.02,0.05) circle (0.035); 
        \fill[gray!40!white] (0.01,-0.06) circle (0.025); 
        \fill[gray!60!white] (-0.06,0.02) circle (0.018);
    \end{scope}

    \node at (0.1,0,0) {\includegraphics[width=0.09\textwidth]{Satellite.png}}; 
    \draw[thick, nicered, ->] (0.1,0,0) -- (0.1,0.2,1.8);  
    \draw[thick, nicered, ->] (0.1,0,0) -- (0.1,-0.6,-0.9); 
    \draw[thick, nicered, ->] (0.1,0,0) -- (0.1,0.4,-0.5); 

\end{tikzpicture}
\end{adjustbox}
\]

Recent work done by Broćić, Cant, and Shelukhin \cite{Broc} shows that, given $R \subset M$ closed manifolds, and $\Lambda \subset \mathbb{S}^\star M$ a Legendrian which is isotopic to $\partial N^\star R$, then \textit{the chord conjecture holds}; meaning that there exists at least one Reeb chord on $(\partial W, \alpha)$ with ends in $\Lambda$. This gives existence of \textit{at least} one trajectory bi-normal to the $x$-axis.

\medskip

However, if one wanted to use the current methods from \cite{Mor01} and \cite{ML24} (namely symplectic/wrapped Floer homologies, or close adaptations like Rabinowitz Floer homology), they would need a way of distinguishing physically relevant chords from the undesirable ones.

\smallskip

Indeed, to define $SH^*$ or $HW^*$, we need to extend our original Hamiltonian $H : W\to\mathbb{R}$ to one on $$\widehat W := W \cup_{\partial W} [1,+\infty)\times\partial W$$

While this can be done fairly easily (see \cite{Mor01}), this process requires making choices, and hence can cause the appearance of undesirable Hamiltonian chords on the extension $[1, +\infty)\times\partial W$.

    \indent For the purposes of Theorem \ref{MLPBThm}, since we were only interested in interior chords, all we needed was a way to ignore chords on the extension $[1,+\infty)\times\partial W$, which was provided by the twist condition $+$ the strong index-definiteness assumption (or, alternatively, the weakened twist condition from \cite{Lim}). Since we are now interested in chords on $\partial W$, we can no longer afford to ignore the extension, and therefore need a way of effectively detecting those chords which are only artifacts of the extension process, and of forcing them to be homologically invisible.
\end{remark}


\bigskip

\begin{thebibliography}{10}
\normalsize

\bibitem[AFvKP12]{AFvKP} Albers, Peter; Frauenfelder, Urs; van Koert, Otto; Paternain, Gabriel P. Contact geometry of the restricted three-body problem. Comm. Pure Appl. Math. 65 (2012), no. 2, 229--263. \url{https://doi.org/10.1002/cpa.21380}.

\smallskip

\bibitem[AFF+12]{Albers2012} Albers, Peter; Fish, Joel W.; Frauenfelder, Urs, et al. Global surfaces of section in the planar restricted 3-body problem. Arch. Rational Mech. Anal. 204 (2012), no. 2, 273--284. \url{https://doi.org/10.1007/s00205-011-0475-2}.

\smallskip

\bibitem[AS08]{ASconormal} Abbondandolo, Alberto; Portaluri, Alessandro; Schwarz, Matthias. The homology of path spaces and Floer homology with conormal boundary conditions. J. Fixed Point Theory Appl. 4 (2008), no. 2, 263--293. \url{https://doi.org/10.1007/s11784-008-0097-y}.

\smallskip

\bibitem[BCS24]{Broc} Broćić, Filip; Cant, Dylan; Shelukhin, Egor. The chord conjecture for conormal bundles. Preprint arXiv:2401.08842, (2024). \url{https://arxiv.org/abs/2401.08842}.

\smallskip

\bibitem[CJK20]{CJK} Cho, WanKi; Jung, Hyojin; Kim, GeonWoo. The contact geometry of the spatial circular restricted 3-body problem. Abh. Math. Semin. Univ. Hambg. 90 (2020), no. 2, 161--181. \url{https://doi.org/10.1007/s12188-020-00222-y}.

\smallskip

\bibitem[Gin10]{Conley_conjecture} Ginzburg, Viktor L. The Conley conjecture. Ann. of Math. (2) 172 (2010), no. 2, 1127--1180. \url{https://doi.org/10.4007/annals.2010.172.1129}.

\smallskip

\bibitem[Lim25]{Lim} Limoge, Arthur. PhD Thesis. Heidelberg University, DE. To appear.

\smallskip

\bibitem[ML24]{ML24} Moreno, Agustin; Limoge, Arthur. A Relative Poincar\'e-Birkhoff theorem. Preprint arXiv:2408.06919, (2024). \url{https://arxiv.org/abs/2408.06919}.

\smallskip

\bibitem[MvK22a]{Mor00} Moreno, Agustin; van Koert, Otto. Global hypersurfaces of section in the spatial restricted three-body problem. Nonlinearity 35 (2022), no. 6, 2920--2970. \url{https://doi.org/10.1088/1361-6544/ac692b}.

\smallskip

\bibitem[MvK22b]{Mor01} Moreno, Agustin; van Koert, Otto. A generalized Poincar\'e-Birkhoff theorem. J. Fixed Point Theory Appl. 24 (2022), no. 2, Paper No. 32, 44 pp. \url{https://doi.org/10.1007/s11784-022-00957-6}.

\smallskip

\bibitem[Mor24]{moreno2024symplecticgeometrythreebodyproblem} Moreno, Agustin. The symplectic geometry of the three-body problem. Preprint arXiv:2101.04438, (2024). \url{https://arxiv.org/abs/2101.04438}.

\smallskip

\bibitem[Poi12]{PLT} Poincar\'e, Henri. Sur un th\'eor\`eme de g\'eom\'etrie. Rend. Circ. Mat. Palermo 33 (1912), 375--407. \url{https://doi.org/10.1007/BF03015314}.


\end{thebibliography}
\end{document}